\documentclass[12pt]{article}

\usepackage[utf8]{inputenc}

\usepackage{amsfonts}
\usepackage{amssymb}
\usepackage{amsmath}
\usepackage{amsthm}
\usepackage[T2A]{fontenc}

\usepackage{amsfonts}
\usepackage{amssymb}
\usepackage{amsmath}
\usepackage{amsthm}

\begin{document}
 
 \begin{Large}
  \centerline{\bf On the distribution of non-primitive}
  \centerline{\bf lattice points in the plane}
  \end{Large}
  \vskip+1cm
  \begin{large}
  \centerline {by {\bf Nikolay Moshchevitin}
\footnote{This work is supported by the Russian Science Foundation under grant 19–11–00001.}}
  \end{large}
    \vskip+1cm
    
    {\bf 1.Introduction.}
    
    \vskip+0.3cm
    
    In \cite{e0} Erd\H{o}s proved two following easy and beautiful theorems.
    
        \vskip+0.3cm
    {\bf Theorem A.}\,
    {\it
    For every $\varepsilon >0$  there exist arbitrarily large positive integer  $x$  and $ y\ge x$ such that 
    $$
    {\rm g.c.d.} (x+i , y+j) >1
    $$
    for all  pairs $i,j$ with
    $$
    0\le i,j \le (1-\varepsilon) \left(\frac{\log x}{\log\log x}\right)^{1/2}.
    $$
    }
        \vskip+0.3cm
    
        {\bf Theorem B.} \,
    {\it  For a certain positive constant $c$ for any positive integers  $ x\le y$ there exist a pair of integers $i,j$ with
    $$
     0\le i,j \le c\,\frac{\log x}{\log\log x}
    $$
    such that 
    $$
    {\rm g.c.d.} (x+i , y+j) =1.
    $$
    }
    
        \vskip+0.3cm
    Applying Theorem B in \cite{e1}
    Chalk and Erd\H{o}s proved the following result on coprime  inhomogeneous approximation to a real number.
        \vskip+0.3cm
    
        {\bf Theorem C.}
        \,
    {\it
    For any given  $\alpha\in \mathbb{R}\setminus\mathbb{Q} $
     and any real number $\eta$, there exists an absolute constant $ C$ such that
 $$
     q\, |q\alpha - \eta - r|  \le C\left( \frac{\log q}{\log\log q}
     \right)^2
     $$
is satisfied by infinitely many coprime integers $ (q,r), q\ge 1$.
    }
    
        \vskip+0.3cm
    
    In the last years the problem on inhomogeneous approximation with coprime numbers again became of interest
    (see \cite{j} and the bibliography therein). In particular, Jitomirskaya and Liu
    proved the following result.
    
        \vskip+0.3cm
    {\bf 
    Theorem D.} 
    \,
  {\it For any constant $C$, there exists 
  $\alpha\in \mathbb{R}\setminus\mathbb{Q} $ and $\eta$ 
  such that the inequality
  $$
q\, | q\alpha - \eta - r|  \le C
  $$
only has finitely many coprime  integer solutions $ (q,r), q\ge 1$.}
    
       \vskip+0.3cm 
    
    In the present paper we improve on the result by Jitomirskaya and  Liu.
    Combining the original argument  of Theorem A by Erd\H{o}s  with the construction from \cite{j}
    we show that  there exist uncountably many $\alpha \in \mathbb{R}\setminus \mathbb{Q}$ and $ \eta\in \mathbb{R}$ such that 
    $$
    \inf_{(q,r)\in \mathbb{Z}^2,\,\,   q> 100, \,\,(q,r) =1}  q \,\,\frac{\log \log q}{\sqrt{ \log q}}\,\,|q\alpha -\eta -r|  >0.
    $$
     In Section  2 below we describe generalized Erd\H{o}s'  construction.
     In Sections 3, 4 we explain how to construct $\alpha$ and $\eta$.
     In Section   5  we give the exact formulation of our result (Theorem 1) and complete the proof.
    
        \vskip+0.3cm
      
    {\bf 2. Constructions with prime numbers.}
    
        \vskip+0.3cm
    Let $ k$ be a non-negative integer. We define the sets of integer points $\Omega_k$ in the following way.
    We put $\Omega_0  = \{ (0,0)\}$ and define 
    $$
    \Omega_k = \{ (i,j)\in \mathbb{Z}^2:\,\, \max (1, |i|)\cdot \max (1, |j|)   \le k\}.
    $$
    Let $\omega_k$ be the cardinality of the set $\Omega_k$.
    It is clear that  $\omega_0 = 1, \omega_1 = 5, \omega_2 = 13$ and 
    by simple counting argument we have  
    \begin{equation}\label{o}
     \omega_k = 4 k\log k + O(k) ,\,\,\,\, k \to \infty .
    \end{equation}

    We arrange the sequence
    $$
    p_1 =2, p_2=3, p_3 = 5, p_4 = 7,...
    $$
    of all prime numbers in the following way.
    Let $ p_{0,0} = 2$ and if 
    we have a permutation of the first $\omega_k$ prime mumbers
    indexed as
    $$
    p_{i,j},\,\,\,\,\,\,\,   (i,j)\in \Omega_k
    $$
    we define next $   \omega_{k+1} - \omega_k$ numbers
    $$
        p_{i,j},
      \,\,\,\,\,\,\,\,  (i,j)\in \Omega_{k+1}\setminus\Omega_k
    $$
    as an arbitrary fixed permutation of prime numbers
    $$
    p_s, \,\,\, \,\,\,\,\omega_k<s \le \omega_{k+1}
.$$
We put 
\begin{equation}\label{p}
P_k =
\prod_{ (i,j)\in \Omega_{k}} p_{i,j} =
\prod_{s=1}^{\omega_k} p_s.
\end{equation}
  By Prime Number Theorem and (\ref{o})  we have
\begin{equation}\label{pp}
{\log P_k}
 \sim  4 k \, (\log k)^2
 ,\,\,\,\,\,
 k \to \infty.
\end{equation}
Define
\begin{equation}\label{Q}
Q_{k+1} = \frac{P_{k+1}}{P_k} ,\,\,\,\,\, {\rm g.c.d.}(P_k, Q_{k+1}) =1.
\end{equation}
For every $ k$ we define numbers
$$
A_i^{[k]} =
\prod_{s: \, (i,s)\in \Omega_k} p_{i,s},\,\,\,\,\, -k \le i \le k
$$
and
$$
B_j^{[k]} =
\prod_{s: \, (s,j)\in \Omega_k}  p_{s,j},
\,\,\,\,\, -k \le j\le k.
$$
It is clear that
$$
 A_i^{[k]} |  A_i^{[k+1]} ,\,\,\,\,\,\,\,
 B_j^{[k]}| B_j^{[k+1]} ,
\,\,\,\,\,\,\,\,
\prod_{i=-k}^k  A_i^{[k]} =
\prod_{j=-k}^k  B_j^{[k]} = P_k,
$$
\begin{equation}\label{ab}
    {\rm g.c.d.}
(A^{[k]}_{i_1}, A^{[k]}_{i_2})=
    {\rm g.c.d.}
(B^{[k]}_{j_1}, B^{[k]}_{j_2}) =1,\,\,\,\,\,
i_1\neq i_2,\,\,\, j_1 \neq j_2
\end{equation}
and 
\begin{equation}\label{ab1}
    {\rm g.c.d.}
(A^{[k_1]}_{i}, B^{[k_2]}_{j}) =
\begin{cases}
p_{i,j},\,\,\, 
\max (1, |i|)\cdot \max (1, |j|) \le \min (k_1,k_2)
\cr
1,\,\,\,\,\,\,\,\,
\max (1, |i|)\cdot \max (1, |j|) > \min (k_1,k_2)
\end{cases}
\end{equation}
for any $k_1, k_2$.

We define $X_k, Y_k$ from the conditions
$$
X_k + i \equiv 0\pmod{A_i^{[k]}},\,\,\, -k \le i \le k,
$$
$$
Y_k + j\equiv 0\pmod{B_j^{[k]}},\,\,\, -k \le j \le k,
$$
$$
P_k\le X_k, Y_k \le 2P_k-1
$$
by Chinese Remainder Theorem. It is clear that
\begin{equation}\label{xy}
    {\rm g.c.d.}
(
X_k + i ,
 Y_k + j
)
\neq 1\,\,\,\,\,
\text{for }\,\,\,\,
\max (1, |i|)\cdot \max (1, |j|) \le k
\end{equation}
and
\begin{equation}\label{xky}
X_{k+1} \equiv X_k \pmod{P_k},\,\,\,\,\,
Y_{k+1} \equiv Y_k \pmod{P_k}.
\end{equation}
We should take into account that 
$$P_0  =X_0 = Y_0=  A_0^{[0]}=  B_0^{[0]}=2.
$$

Also  for non-collinear  integer vectors $\pmb{f},\pmb{g}\in \mathbb{Z}^2$ we need to consider  affine lattices
$
\Lambda_{k} (
\pmb{f},\pmb{g}
)
$ defined for even  values of $k$ as
$$
\Lambda_{2l} (
\pmb{f},\pmb{g}
)
=
\{ \pmb{z} = x\pmb{f}+y\pmb{g} :\,\,\,
x\equiv X_l\pmod{P_l},\,\,\, 
y\equiv Y_l\pmod{P_l}\}.
$$
and for odd
values of $k$ as
$$
\Lambda_{2l+1} (
\pmb{f},\pmb{g}
)
=
\{ \pmb{z} = x\pmb{f}+y\pmb{g} :\,\,\,
x\equiv Y_{l+1}\pmod{P_{l+1}},\,\,\, 
y\equiv X_l\pmod{P_l}\}.
$$
    \newpage
  {\bf 3. Continued fractions and integer points.}
  
      \vskip+0.3cm
Now we consider  irrational $\alpha$ defined by its continued fraction expnsion
\begin{equation}\label{cf}
\alpha = [a_0;a_1,a_2,a_3,...,a_k,...],
\,\,\,\,\, a_0 \in \mathbb{Z},\,\,\, a_j \in \mathbb{Z}_+, \, j = 1,2,3,...  .
\end{equation}
We suppose in our construction that
\begin{equation}
\label{pe}
a_{k+1} \equiv 0 \pmod{P_{[k/2]}} 
\end{equation}
and
\begin{equation}
\label{boo}
a_{k+1} \ge  W_k \cdot {{P}_{[k/2]+1}} 
\,\,\,\,\,\,
\text{with integer}\,\,\,\,\,\, W_k \to \infty,\,\, k\to \infty.
\end{equation}
Let 
$$\frac{u_k}{v_k} 
 = [a_0;a_1,a_2,a_3,...,a_k],\,\,\,\,     {\rm g.c.d.} (v_k,u_k) = 1
$$ 
be convergent fractions to $\alpha$.
Define  primitive vectors 
$$ \pmb{e}_k = (v_k, u_k)\in \mathbb{Z}^2,\,\,\, k = 0,1,2,3,... .
$$
It is clear that 
\begin{equation}\label{qqq}
\pmb{e}_{k+1} = a_{k+1} \pmb{e}_k + \pmb{e}_{k-1}.  
\end{equation}
It would be
convenient to define the lattice
$$
\frak{L}_k = \Lambda_{k} (\pmb{e}_{k-1},\pmb{e}_k)
,
$${
so for $ l = [k/2]$ we have }
$$
\frak{L}_{2l} =  
\{ \pmb{z} = x\pmb{e}_{2l-1}+y\pmb{e}_{2l} :\,\,\,
x\equiv X_l\pmod{P_l},\,\,\, 
y\equiv Y_l\pmod{P_l}\},
 $$
 $$
\frak{L}_{2l+1} =  
\{ \pmb{z} = x\pmb{e}_{2l}+y\pmb{e}_{2l+1} :\,\,\,
x\equiv Y_{l+1}\pmod{P_{l+1}},\,\,\, 
y\equiv X_l\pmod{P_l}\}.
$$

Now we define inductively integer points
\begin{equation}\label{ze}
\frak{Z}_k = (b_k,c_k) = \frak{X}_k \pmb{e}_{k-1}+ \frak{Y}_k \pmb{e}_k \in \frak{L}_k, \,\,\,
k=0,1,2,3,... .
\end{equation}
Put
$
\frak{Z}_0  = (0,0)
$. Then, if  $
\frak{Z}_k \in \frak{L}_k
$
is defined, we define 
$
\frak{Z}_{k+1} \in \frak{L}_{k+1}
$
by the following procedure.
We take $ l =[k/2]$ and define
\begin{equation}\label{z0}
\frak{Z}_{k+1} =
\frak{Z}_k
+
 {P}_{l}\lambda_{k+1} \pmb{e}_k =
\frak{X}_k \pmb{e}_{k-1}
+
(
\frak{Y}_k+
 {P}_{l}\lambda_{k+1} )\pmb{e}_k  
 ,
\end{equation}
so
\begin{equation}\label{z2}
{b}_{k+1} = {b}_{k}+ {P}_{l}\lambda_{k+1}v_k\,\,\,\,\,
\text{and}
\,\,\,\,\,
{c}_{k+1}  =  {c}_{k}+ {P}_{l} \lambda_{k+1}u_k.
\end{equation}
The value of $\lambda_{k+1} \in \mathbb{Z}_+$  will be defined later in (\ref{bound}).
 
Taking into account equality 
(\ref{qqq}), by (\ref{ze},\ref{z0}) we can write
\begin{equation}\label{z11}
\frak{Y}_{k+1} = \frak{X}_{k}
\equiv
\begin{cases}
{X}_{l} \pmod{P_l} , \,\,\,\,\,\,\,\,\,\text{if $k$ even}
\cr
Y_{l+1} \pmod{P_{l+1}}  , \,\,\text{if $k$ odd}
\end{cases}.
\end{equation}
{and}
\begin{equation}\label{z12}
\frak{X}_{k+1}  =  \frak{Y}_{k}+{P}_{l}\lambda_{k+1} -\frak{X}_k a_{k+1},
\end{equation}

By (\ref{pe}) we have $\frak{X}_{k+1} \equiv \frak{Y}_{k}\pmod{P_l}$.
Since $
\frak{Z}_k \in \frak{L}_k$, by (\ref{xky}) we see that 
\begin{equation}\label{raz}
\frak{X}_{k+1} \equiv 
\begin{cases}
{Y}_{l} \equiv  {Y}_{l+1} , \,\,\text{if $k$ even}
\cr
X_{l} \equiv X_{l+1} , \,\,\text{if $k$ odd}
\end{cases}
\pmod{P_l}.
\end{equation}
On the other hand  $    {\rm g.c.d.} (P_l, Q_{l+1}) = 1$ (see (\ref{Q})). So we can find $ \lambda^*_{k+1} \in \{1,2,...,Q_{l+1}\}$ such that for every $\lambda_{k+1} \equiv \lambda^*_{k+1} \pmod{Q_{l+1} }$  one has
\begin{equation}\label{dva}
\frak{X}_{k+1} \equiv 
\begin{cases}
 {Y}_{l+1} , \,\,\text{if $k$ even}
\cr
X_{l+1} , \,\,\text{if $k$ odd}
\end{cases}
\pmod{Q_{l+1}}.
\end{equation}
As
$P_{l+1} = P_{l} Q_{l+1}$
equalities (\ref{raz},\ref{dva}) imply
\begin{equation}\label{tri}
\frak{X}_{k+1} \equiv 
\begin{cases}
 {Y}_{l+1} , \,\,\text{if $k$ even}
\cr
X_{l+1} , \,\,\text{if $k$ odd}
\end{cases}
\pmod{P_{l+1}}.
\end{equation}
We see that equalities (\ref{z11}) and (\ref{tri}) ensure   $\frak{Z}_{k+1}\in \frak{L}_{k+1}$.

  Finally we take $\lambda_{k+1} $ of the form  
  $\lambda_{k+1} =  \lambda_{k+1}^*+\lambda Q_{l+1}, \lambda 
  \in \mathbb{Z}$ 
  to satisfy 
\begin{equation}\label{bound}
\left|
 \lambda_{k+1} - \frac{a_{k+1}}{2{P}_{l}} \right| \le
 \frac{Q_{l+1}}{2} ,
\end{equation}
so $\lambda$ is the nearest integer to $\frac{a_{k+1}}{2{P}_{l+1}}-\frac{ \lambda_{k+1}^*}{Q_{l+1}}$.

We have defined $\frak{Z}_{k+1} = (b_{k+1}, c_{k+1}) \in \frak{L}_{k+1}$ of the form 
(\ref{z0}).

 Here we should not that from (\ref{z2}) and (\ref{bound}) it follows that
 $$
 \left|
 b_{k+1}-\left( b_k +\frac{v_{k+1}}{2}\right)\right| \le P_{l+1} v_k
 $$
 and by induction we have
 \begin{equation}\label{be}
 \left| b_{k+1} - \frac{v_{k+1}}{2}\right|\le 2P_{l+1}v_k.
 \end{equation}
   
       \vskip+0.3cm
{\bf 4. The value of $\eta$.}
    \vskip+0.3cm

From (\ref{z2}) we see that 
$$ b_{k} = \sum_{j=1}^k  P_{[(j-1)/2]}\lambda_j  v_{j-1},\,\,\,\,\,
c_{k} = \sum_{j=1}^k  P_{[(j-1)/2]}\lambda_j  u_{j-1}.
$$
Now we define
\begin{equation}\label{eta}
\eta = \lim_{k\to\infty} ( b_k\alpha-c_k) =
\sum _{j=1}^\infty P_{[(j-1)/2]} \lambda_j ( v_{j-1}\alpha-u_{j-1}).
\end{equation}
Of course, we should explain why the series in  (\ref{eta}) converges.
We use the estimate
$$
 \frac{1}{v_j+v_{j-1}}\le
|
  v_{j-1}\alpha - u_{j-1} |  \le \frac{1}{v_j},
$$
By means of the upper bound here  we see that 
$$
|
 P_{[(j-1)/2]+1}
\lambda_j ( v_{j-1}\alpha-u_{j-1} )|\le
\frac{  P_{[(j-1)/2]} \lambda_j}{v_j} 
\le
\frac{  1}{ 2v_{j-1}} + \frac{ P_{[(j-1)/2]}}{ 2v_{j}}\le
\left( \frac{1}{2} + \frac{1}{2W_j}\right)\,\frac{1}{v_{j-1}}
,
$$
by (\ref{bound}) and (\ref{boo},\ref{qqq}), as
$\frac{v_{j}}{v_{j-1}}\ge a_j \ge W_j P_{[(j-1)/2]}$. So the series from (\ref{eta}) converges absolutely since
$ v_j \ge \left(\frac{1+\sqrt{5}}{2}\right)^{j-1}$. 
Moreover, it is easy to see that 
the absolute values of the summands from series (\ref{eta}) decrease monotonically to zero.
Futhermore, as $ P_{[(j-1)/2]} \lambda_j $ is positive and
$u_{j-1} - v_{j-1}\alpha$ changes its sign we see that 
\begin{equation}\label{eta1}
|\eta -   ( b_k\alpha-c_k) |<
\left( \frac{1}{2} + \frac{1}{2W_k}\right)\,\frac{1}{v_{k}}
 .
\end{equation}
The lower bound for the approximation is given by
$$
|\eta -   ( b_k\alpha-c_k) |>
 P_{[k/2]} \lambda_{k+1}|u_{k} - v_{k}\alpha|-
 P_{[(k+1)/2]} \lambda_{k+2}|u_{k+1} - v_{k+1}\alpha|
\ge
$$
\begin{equation}\label{eta2}
\ge
\frac{P_{[k/2]} \lambda_{k+1}}{v_{k+1}+v_k}
-
\frac{P_{[(k+1)/2]} \lambda_{k+2}}{v_{k+2}}
\ge 
\left( \frac{1}{2} - \frac{2}{W_k}\right)\,\frac{1}{v_{k}} 
\end{equation}
(for the second summand we used here (\ref{eta1})).
In addition by  (\ref{qqq})  and (\ref{boo}) we see that 
$$
|\eta -  ( b_k\alpha-c_k) |>|\eta - ( b_{k+1}\alpha-c_{k+1}) |.
$$
We should note that the differences
$$
\eta -  ( b_k\alpha-c_k)
\,\,\,\,\,
\text{and}
\,\,\,\,\,
\eta -    ( b_{k+1}\alpha-c_{k+1}) 
$$
have different signs meanwhile the differences
$$
\eta -    ( b_k\alpha-c_k)
\,\,\,\,\,
\text{and}
\,\,\,\,\,
\eta -   ((b_{k+1}-v_k)  - (c_{k+1}-u_k)\alpha)
$$
have the same sign.

     \vskip+0.3cm
 
 {\bf 5. The main result.}
     \vskip+0.3cm
 
 Now we are able to give the exact formulation of our main result.
 
     \vskip+0.3cm
 {\bf Theorem 1.}
 {\it Suppose that the partial quotients of the continued fraction expansion (\ref{cf}) for $\alpha$
 satisfy (\ref{pe}) and (\ref{boo}). 
 Suppose that with some positive $\gamma$ the inequality
 \begin{equation}\label{iiee}
 a_k \le \exp (\gamma k (\log k)^2)
 \end{equation}
 holds.
 Let $\eta$ be defined by (\ref{eta}).
 Then for any primitive point $(q,r)\in \mathbb{Z}^2$ with  $     {\rm g.c.d.} (q,r) =1$ and $q >100$ one has
 \begin{equation}\label{main}
 q\, |q\alpha -\eta-r| > c \,\frac{\sqrt{\log q}}{\log \log q},
 \end{equation}
 with certain positive $c$.
 }
 
     \vskip+0.3cm
  The proof of Theorem 1 is close to original argument from \cite{j}, as well as to the argument  by Worley \cite{w}
 (see also \cite{d}).
 
 We suppose that constants in symbols $\gg$ and $ \asymp$ below depend on $\gamma$.
 
     \vskip+0.3cm
 {\bf Remark 1.} In is clear that for any $\alpha$  under conditions  (\ref{pe},\ref{boo},\ref{iiee})
 it is possible to find uncountably many $\eta$ satisfying the conclusion of Theorem 1.
     \vskip+0.3cm
 
  {\bf Remark 2.}   We can choose $a_k$ satisfying (\ref{pe},\ref{boo}) and  (\ref{iiee})
  because of  (\ref{pp}).
  From the inequalities
  $
  v_k \le 2^k \prod_{j=1}^ k a_k
  $
  and (\ref{iiee}) we deduce that 
  $$
   k \gg  \frac{\sqrt{\log v_k}}{\log\log v_k}.
   $$
   Moreover, from (\ref{be}) it follows that $ b_k \asymp v_k$ and so by (\ref{boo}) and (\ref{iiee}) we have
  \begin{equation}\label{rr}
 k \gg  \frac{\sqrt{\log b_k}}{\log\log b_k} \asymp  \frac{\sqrt{\log b_{k+1}}}{\log\log b_{k+1}}.
 \end{equation}
  
      \vskip+0.3cm

 Proof of Theorem 1.
 
 From (\ref{rr}) we see that 
  it is enough to prove the inequality 
   \begin{equation}\label{rrr}
    q\,|q\alpha - \eta - r| \gg
 k
 \end{equation}
  for every primitive point 
  $ \pmb{z} =(q,r) \in \mathbb{Z}^2$ 
  satisfying
  \begin{equation}\label{kaa}
  b_k \le q< b_{k+1}.
  \end{equation}

  Each point $ \pmb{z} =(q,r) \in \mathbb{Z}^2$ 
 can be written
   in the form
  \begin{equation}\label{kaa1}
   \pmb{z}=  \frak{Z}_k+ x\pmb{e}_{k-1}+ y\pmb{e}_k
   \end{equation}
   and in another form
    \begin{equation}\label{kaa2}
    \pmb{z}=  \frak{Z}_{k+1}+ x'\pmb{e}_{k}+ y'\pmb{e}_{k+1}
   \end{equation}
    with integer $x,y,x',y'$.
  
 We should note that  if an integer point  $ \pmb{z} =(q,r) $ can be represented in the form
 \begin{equation}\label{from}
 \pmb{z}=  \frak{Z}_k+ x\pmb{e}_{k-1}+ y\pmb{e}_k,\,\,\,\,
 \text{with}\,\,\,\, \max (1, |x|)\cdot \max (1, |y|) \le \frac{k}{2},
 \end{equation}
 then 
 by (\ref{xy}) and by the condition $\frak{Z}_k \in \frak{L}_k$ we see  that  for all $\pmb{z}$ of the form 
 (\ref{from}) we have $ (\frak{X}_k+x,\frak{Y}_k+y)\neq 1$ and so $ (q,r)\neq 1$.
 
 Now we consider integer points $ \pmb{z} =(q,r) $  under the condition (\ref{kaa}) and such that simultaneously 
 \begin{equation}\label{or}
\max (1, |x|)\cdot \max (1, |y|)  > \frac{k}{2}\,\,\,\,\text{in (\ref{kaa1})},
 \,\,\,\,
 \text{and}
 \,\,\,\,
\max (1, |x'|)\cdot \max (1, |y'|)  > \frac{k}{2}\,\,\,\,\text{in (\ref{kaa2})}.
 \end{equation}
 For these points we shall prove that  (\ref{rrr}) holds and this will complete the proof of Theorem 1.
 
 Put $l =\left[\frac{k}{2}\right]$ and consider two cases:

  \noindent 
 {\bf case $1^0$}: $ \pmb{z} = (q,r) = \frak{Z}_k +y \pmb{e}_k$ with $ l \le y \le P_{k}\lambda_{k+1} - l$;

 \noindent and

  \noindent 
  {\bf case $2^0$}:  $ \pmb{z} = (q,r) = \frak{Z}_k +  x\pmb{e}_{k-1} +y \pmb{e}_k$  with 
  $  x\neq 0$.
  
  As $ \frak{Z}_{k+1} = (b_{k+1},c_{k+1}) = \frak{Z}_k  + P_{k}\lambda_{k+1}\pmb{e}_k$,
  in   {\bf case $1^0$} inequality  (\ref{kaa}) is satisfied automatically.
In   {\bf case $2^0$} we  should assume (\ref{kaa}).

   In   {\bf case $1^0$},
 $$
   q = b_k + y v_k, \,\,\,\,\,\,
   ||\eta - q\alpha|| =   \eta -   (  (b_k+yv_k)\alpha - (c_k  +yu_k) )
  $$
   and
   the differences 
   $$
   \eta -   (   (b_k+lv_k)\alpha - (c_k  +lu_k))
\,\,\,\,\,
\text{and}
\,\,\,\,\,
\eta -   (  (b_{k+1}-lv_k)\alpha - (c_k  +lu_k))
   $$
   have the same sign (see the very end of Section  4). 
   
   We are interested in getting lowed bound for approximation for primitive points $\pmb{z}$, and 
   it is enough to consider only those $\pmb{z}$ for which (\ref{or}) holds. So we can assume that 
   $$
   b_k+lv_k\le q\le  b_{k+1}-lv_k.
   $$
   For
   $$
   \pmb{z} = (q,r) = (b_k+lv_k, c_k + lu_k)
   $$
   we have $ q = b_k+lv_k,> lv_k$ and
   $$
   q  |\eta -   q\alpha - r|\ge 
   lv_k (|\eta - (b_k\alpha+c_k)| - l |v_k\alpha - u_k|)\ge
   \frac{l}{5}-   \frac{l^2 v_k}{v_{k+1} }\ge   \frac{l}{10}-   \frac{l^2 }{a_{k+1} }\ge    \frac{l}{20}
   $$
   as $
|
  v_{k}\alpha - u_{k} |  \le \frac{1}{v_{k+1}}
$
and
   by (\ref{eta2})  and (\ref{boo}).
   
   For
   $$
   \pmb{z} = (q,r) = (b_{k+1}-lv_k, c_{k+1}-lu_k)
   $$
   we have $ q= b_{k+1}-lv_k\ge b_{k+1}/2\ge v_{k+1}/8$ (we use (\ref{be})) and
   $$
   q  |\eta -   q\alpha - r|\ge 
   \frac{v_{k+1}}{8}
 (  l |v_k\alpha - u_k|- |\eta - (b_{k+1}\alpha+c_{k+1})|)\ge
   \frac{l}{10}-   \frac{1}{8 } \ge    \frac{l}{20},
   $$ 
   for $ l\ge 3$,
   as  we have inequalities $
|
  v_{k}\alpha - u_{k} |  \ge \frac{1}{v_{k+1}+v_k}
$
   and (\ref{eta1}).

Now
   we see that for $\pmb{z} = (q,r)=  (c_k  +yu_k, b_k+yv_k)$   the inequality
   $$
   q\,
     |\eta - q\alpha-r| =
    (b_k + y v_k) |\eta -   ((c_k  +yu_k) - (b_k+yv_k)\alpha)  |\ge    \frac{l}{20}
   $$
   is true for any intermediate value of $y$ satisfying $ l \le y \le P_{k}\lambda_{k+1} - l$  by convexity argument.
    Now we take into account (\ref{rr}) and (\ref{kaa}).
   In {\bf case $1^0$} for primitive $\pmb{z}$ inequality (\ref{rrr}) is proven.
   
        The proof in {\bf case $2^0$}  is based on the same ideas. In this case we have $ x\neq 0$ and for any $\pmb{z} = (q,r)$   under  consideration from   the lower bound (\ref{eta2}) we deduce the inequality
  \begin{equation}\label{c0}
         |\eta - q\alpha-r| 
     \ge \frac{|x|}{3v_k}. 
     \end{equation}

     By (\ref{be}) and (\ref{kaa})  
     \begin{equation}\label{az}
      q\ge b_k\ge v_k/4.
      \end {equation}
      
       If $ |x|\ge l$, then  by (\ref{az}) we see that 
      $$
 q 
        \, |\eta - q\alpha-r|  \ge  \frac{l}{12}. 
        $$

 Suppose $ 0\neq |x|< l$.   
 As we deal with primitive $\pmb{z}$     we should not consider $\pmb{z}$ which   satisfy (\ref{from}).  So 
 $ y>\frac{l}{|x|}
 $ 
 and 
 instead of (\ref{az}) we have
             \begin{equation}\label{c1}
              q\ge   yv_{k} - b_k - l v_{k-1} \ge \frac{yv_{k}}{4} >\frac{lv_{k}}{4|x|}
              .
              \end{equation}
              Now inequalities (\ref{c0},\ref{c1}) give us
              $$
                q\,
     |\eta - q\alpha-r| \ge  \frac{l}{12}.
     $$
     
     In {\bf case 2}  inequality (\ref{rrr}) is proven also. 
     
     Theorem 1 is proven.
     
         \vskip+0.3cm
         The author thanks Oleg German and Vasilii Neckrasov for  fruitful discussions.

            \vskip+0.3cm
            author:
            
            Nikolay G. Moshchevitin
            
         Steklov Mathematical Institute, 
         
         ul. Gubkina 8, Moscow,  119991
         
         Russia

\end{document}